\documentclass{bmcart}

\usepackage[utf8]{inputenc} 
\usepackage{amssymb}
\usepackage{amsthm}

\def\includegraphics{}

\usepackage{graphics}
\usepackage{epsfig}
\usepackage{graphicx}
\usepackage{epstopdf}

\usepackage{amsmath, amsthm, amscd, amsfonts, amssymb, graphicx, color}
\usepackage{amsthm}
\usepackage{amssymb}
\usepackage{multirow,multicol}
\usepackage{booktabs}
\usepackage{amsmath, amscd, amsfonts, amssymb, graphicx, color}

\usepackage{float}

\begin{document}

\begin{frontmatter}

\begin{fmbox}
\dochead{Research}


\title{Symplectic method for Hamiltonian stochastic differential equations with multiplicative L\'{e}vy  noise in the sense of Marcus}


\author[
   email={zhan2017@fafu.edu.cn}   
]{\fnm{Qingyi} \snm{Zhan}}$^{1}$
\author[
]{Jinqiao Duan$^2$}
\author[
]{Xiaofan Li$^2$}
\author[
]{Yuhong Li$^3$}



\begin{artnotes}

\end{artnotes}

\end{fmbox}


\begin{abstractbox}

\begin{abstract}
  A class of Hamiltonian stochastic differential equations with multiplicative L\'{e}vy  noise in the sense of Marcus, and the construction and numerical implementation methods of symplectic Euler scheme, are considered. A general symplectic Euler scheme for this kind of Hamiltonian stochastic differential equations is devised, and its convergence theorem  is proved. The second part presents realizable numerical implementation methods for this scheme in details. Some numerical experiments are conducted to demonstrate the effectiveness and superiority of the proposed method by the simulations of its orbits, Hamlitonian,and convergence order over a long time interval.

\end{abstract}




\end{abstractbox}
%

\end{frontmatter}



\section{Introduction}

Recently, there have been increasing interests in the stochastic differential equations(SDEs) with non-Gaussian noise. These SDEs have played an important role in the theory and application of stochastic dynamics\cite{D. Applebaum, J.Duan}. In the research of some Hamiltonian SDEs with L\'{e}vy noise in the sense of Marcus form, which can preserve the symplectic structure, many people have paid more and more attention in numerical simulations of random  phenomena\cite{K. Feng, E.Hairer,J.Hong}. 
Numerical computations are crucial to study dynamical behaviour of Hamiltonian SDEs. And a theoretical framework, structure-preserving algorithm, has been widely applied in many aspects.
 Therefore, we investigate the reliability  and feasibility of numerical computations of Hamiltonian SDEs with multiplicative L\'{e}vy noise.

 This work deals with symplectic Euler scheme of a class of Hamiltonian SDEs with multiplicative L\'{e}vy  noise in the sense of Marcus. This is a continuation of \cite{Q.Zhan6}, where symplectic Euler scheme of the Hamiltonian SDEs with additive L\'{e}vy  noise in the sense of Marcus was considered. It is motivated by two facts. Firstly, as we know, in the case of deterministic Hamiltonian differential equations, symplectic methods are available\cite{K. Feng, E.Hairer}. Therefore, it is natural to expect to construct symplectic methods for Hamiltonian SDEs with multiplicative L\'{e}vy  noise, which are also needed in many aspects. The related works about Hamiltonian SDEs with Gaussian noise are shown in \cite{G.Milstein}. Many contributions are made to the numerical analysis of SDEs \cite{T. Wang, X. Wang}, and  numerical methods of SDEs can be seen in \cite{Milstein},\cite{Q.Zhan}-\cite{Q.Zhan5}. Secondly,
 the construction of conditions which can preserve the Hamiltonian structure of SDEs driven by additive L\'{e}vy  noise has been finished in \cite{Q.Zhan6}.  These results are the foundations of symplectic scheme of Hamiltonian SDEs with multiplicative  L\'{e}vy  noise.
 To the best of our knowledge, no investigations of symplectic scheme of Hamiltonian SDEs with multiplicative  L\'{e}vy  noise in the sense of Marcus  exist in the literature so far.

In this work, we focus on symplectic Euler scheme for Hamiltonian SDEs with multiplicative L\'{e}vy noise in the sense of Marcus form. We compare the numerical dynamical behaviors of symplectic Euler method with non-symplectic methods in these aspects, which include the Hamiltonian, the preservations of symplectic structure and the orbits in a long time interval. All these are shown in the numerical experiments. For our purpose that the numerical experiments are realizable and  simply achieved by programming, the L\'{e}vy noise is still restricted to be compound Poisson noise with a special realization\cite{T. Li}.

 Our results show that under certain appropriate assumptions the Hamiltonian is approximatively preserved in the discrete time case due to the discontinuous input of multiplicative L\'{e}vy noise. The numerical solution by this scheme can simulate the dynamical behaviour of Hamiltonian SDEs more accurately than non-symplectic methods in a long time interval.

The paper is organized as follows. Section 2 deals with some preliminaries. In Section 3 the theoretical results of preservation of symplectic structure are summarized. The mean-square convergence theorem of this scheme is proved in Section 4. Section 5 presents the details of the numerical implementations for the symplectic Euler scheme. Illustrative numerical experiments are included in Section 6. Finally, the last section is addressed to summarize the conclusions of the paper.

\section{Preliminaries}

We consider the following Hamiltonian SDE with multiplicative non-Gauss L\'{e}vy noises in the sense of Marcus on $\mathbb{M}$,
\begin{equation}\label{2.1}
 dX(t)=V_0(X(t))dt+\sum_{r=1}^m V_r(X(t)) \diamond  dL^r(t),\ \ \  \ \ \ X(0):=X(t_0)=x,
\end{equation}
where $X\in \mathbb{R}^d$,$V_r: \mathbb{R}^d \rightarrow\mathbb{R}^d,r=0,1,...,m $,is the Hamiltonian vector fields, and $\mathbb{M}$ is a smooth $d$-dimensional manifold. L\'{e}vy noises $L$ is defined on a filtered probability space $(\Omega, \mathcal{F}, \{\mathcal{F}_t\}_{t\geq 0}, \mathbb{P})$.
Let $L(t)$ be a $d$-dimensional L\'{e}vy process with the generating triplet $(\gamma,A,\nu)$,
where $\gamma$ is a $d$-dimensional drift vector, $A$ is a symmetric non-negative definite $d \times d$ matrix, and $\nu $ is a radially symmetric L$\acute{e}$vy jump measure on $\mathbb{R}^d \backslash 0$.
Here the Marcus integral for SDE(\ref{2.1}) through Marcus mapping is usually written as
$$X(t)=x+\int_0^t V_0(X(s))ds+\sum_{r=1}^m \int_0^t V_r(X(s-)) \diamond  dL^r(s),\ \ \  \ \ \ $$
 which is defined as
$$X(t)=x+\int_0^tV_0(X(s))ds$$
$$+\sum_{r=1}^m\int_0^tV_r(X(s-))\circ dL_c^r(s) +\sum_{r=1}^m\int_0^tV_r(X(s-))dL_d^r(s)  $$
$$+\sum_{r=1}^{m}\sum_{0\leq s\leq t}[\Phi^r(\Delta L^r(s), V_r(X(s-)),X(s-))-X(s-)-V_r(X(s-))\Delta L^r(s)],$$
where $L_c(t)$ and $L_d(t)$ are the usual continuous and discontinuous parts of $L(t)$, that is, $L(t)=L_c(t)+L_d(t)$. The notation $\circ$  denotes the Stratonovitch differential. And the flow map $\Phi^r(l,v(x),x)$ is the value at $s=1$ of the solution defined through the ordinary differential equations
\begin{equation}
 \left\{
  \begin{array}{lcr}
   \frac{d \xi^r}{ds}=V_r(\xi^r)l, s\in[0,1],\\
   \\
    \xi^r(0)=x.
  \end{array} \right.
 \end{equation}

Let us write Hamiltonian SDEs of even dimension $d=2n$ in the form of
\begin{equation}\label{2.2}
\begin{split}
dP=-\frac{\partial H_0}{\partial Q}(P,Q)dt-\sum_{r=1}^m \frac{\partial H_r}{\partial Q}(P,Q)\diamond  dL^r(t),\ \ \  \ \ \ P(t_0)=p,\\
dQ=\frac{\partial H_0}{\partial P}(P,Q)dt+\sum_{r=1}^m \frac{\partial H_r}{\partial P}(P,Q)\diamond  dL^r(t),\ \ \  \ \ \ Q(t_0)=q,
 \end{split}
\end{equation}
where $X=(P,Q)$, $X_0=(p,q)$  and
$V_r=(\frac{\partial H_r}{\partial P},-\frac{\partial H_r}{\partial Q})$, $r=0,1,2,...,m.$ Here, $P,Q,p,q$ are $n$-dimensional column-vectors.
 We assume that the functions $V_r,r=0,1,2,...,m $ satisfy the conditions which is the same as the conditions in \cite{P. Wei} such that Hamilton SDEs $(\ref{2.2})$
 have an unique global solution, and the solution process is adapted and c$\grave{a}$dl$\grave{a}$g.

We introduce the following notations.

Let $\mathbb{L}^2(\Omega,\mathbb{P})$ be the space of all bounded square-integrable random variables $x:\Omega \rightarrow \mathbb{R}^d$. For random vector $x=({x_1},{ x_2},...,{ x_d}) \in \mathbb{R}^d$, the norm of $x$ is defined in the form of
\begin{equation}\label{2.3}
\|x\|_2=\Big[\int_\Omega[|x_1(\omega)|^2+|x_2(\omega)|^2+...+|x_d(\omega)|^2]d\mathbb{P}\Big]^{\frac{1}{2}}< \infty.
\end{equation}
We define the norm of random matrices as follows
\begin{equation}\label{2.4}
 \| G \|_{\mathbb{L}^2(\Omega,\mathbb{P})} =\Big[\mathbb{E}(|G|^2)\Big]^{\frac{1}{2}},
\end{equation}
where $G$ is a random matrix and $|\cdot|$ is the operator norm.

For simplicity, the norms $\|\cdot\|_2$ and $ \| \cdot \|_{\mathbb{L}^2(\Omega,\mathbb{P})}$ are usually written as $\|\cdot\|$.

\section{Theoretical results on preservation of symplectic structure }

\subsection{Preservation of symplectic structure for continuous Hamiltonian SDEs}
The following lemma is from \cite{P. Wei}.

{\bf Lemma 3.1}\emph{
The Hamiltonian SDE $(\ref{2.2})$ preserves symplectic structure, that is,
$$dP\wedge dQ=dp\wedge dq, i.e., \sum_{i=1}^ndP^i\wedge dQ^i= \sum_{i=1}^ndp^i\wedge dq^i,$$
where $dP\wedge dQ$ is a differential two-form,$P=(P^1,P^2,...,P^n),$ and $Q=(Q^1,Q^2,...,Q^n)$.}

This lemma shows the  preservation of symplectic structure for Hamiltonian SDEs in the case of continuous time.

\subsection{Preservation of symplectic structure for discrete Hamiltonian SDEs}
In this section we consider the Hamiltonian SDE with multiplicative L\'{e}vy  noise as follows,
\begin{equation}\label{3.1}
\begin{split}
 dP=-\sigma_0(P,Q)dt-\sum_{r=1}^m \sigma_r(P,Q)\diamond  dL^r(t),\ \ \  \ \ \ P(t_0)=p,\\
 dQ=\gamma_0(P,Q)dt+\sum_{r=1}^m \gamma_r(P,Q)\diamond  dL^r(t),\ \ \  \ \ \ Q(t_0)=q,
 \end{split}
\end{equation}
where $$\sigma_0(P,Q)= \frac{\partial H_0}{\partial Q}(P,Q),\gamma_0(P,Q)=\frac{\partial H_0}{\partial P}(P,Q),$$
$$\sigma_r(P,Q)=\frac{\partial H_r}{\partial Q}(P,Q),\gamma_r(P,Q)=\frac{\partial H_r}{\partial P}(P,Q),r=1,2,...,m.$$

We make the following assumptions.

{\bf Assumption 1.}

\emph{The drift functions $\sigma_0$ and $\gamma_0$ satisfy the Lipschitz condition
$$|\sigma_r(X_1)-\sigma_r(X_2)|\leq K|X_1-X_2|, |\gamma_r(X_1)-\gamma_r(X_2)|\leq K|X_1-X_2|,$$
where $K$ is a constant,  and $X_i=(P_i,Q_i)\in \mathbb{M},i=1,2, r=0,1,...,m$.}

The exact solution $X_{t_j}:=(P_{t_j},Q_{t_j})$ of $(\ref{3.1})$ at the time $t_j$ is shown as
 \begin{equation}\label{3.3}
\begin{split}
 P_{t_{j+1}}=P_{t_j}- \int_{t_j}^{t_{j+1}} \sigma_0(P(s),Q(s))ds-\sum_{r=1}^m \int_{t_j}^{t_{j+1}} \sigma_r(P(s),Q(s)) \diamond  dL^r(s),\ P(t_0)=p,\\
 Q_{t_{j+1}}=Q_{t_j}+\int_{t_j}^{t_{j+1}} \gamma_0(P(s),Q(s))ds+\sum_{r=1}^m \int_{t_j}^{t_{j+1}} \gamma_r(P(s),Q(s)) \diamond  dL^r(s),\ Q(t_0)=q,
 \end{split}
\end{equation}
where the Marcus integral for SDE(\ref{3.3}) is usually defined by
$$\int_{t_j}^{t_{j+1}} \sigma_r(P(s),Q(s)) \diamond  dL^r(s)$$
$$=\int_{t_j}^{t_{j+1}}\sigma_r(P(s),Q(s))\circ \Delta L_c^r(s)+\int_{t_j}^{t_{j+1}}\sigma_r(P(s),Q(s))\Delta L_d^r(s)$$
$$+\sum_{t_j\leq s\leq t_{j+1}}\Big[\Phi_1^r(\Delta L^r(s),\sigma_r(P(s),Q(s)),P(s-))-P(s-)-\sigma_r(P(s),Q(s))\Delta L^r(s)\Big],$$
and
$$\int_{t_j}^{t_{j+1}} \gamma_r(P(s),Q(s)) \diamond  dL^r(s)$$
$$=\int_{t_j}^{t_{j+1}}\gamma_r(P(s),Q(s))\circ \Delta L_c^r(s)+\int_{t_j}^{t_{j+1}}\gamma_r(P(s),Q(s))\Delta L_d^r(s)$$
$$+\sum_{t_j\leq s\leq t_{j+1}}\Big[\Phi_2^r(\Delta L^r(s),\gamma_r(P(s),Q(s)),Q(s-))-Q(s-)-\gamma_r(P(s),Q(s))\Delta L^r(s)\Big].$$

And the flow maps $\Phi_1^r(l,\sigma_r(P(s),Q(s)),P(s-))$ and $\Phi_2^r(l,\gamma_r(P(s),Q(s)),Q(s-))$ are the value at $\hat{s}=1$ of the solutions defined through the ordinary differential equations, respectively,
\begin{equation}
 \left\{
  \begin{array}{lcr}
   \frac{d \xi_1^r}{d\hat{s}}=\sigma_r(\xi_1^r)l, \xi_1^r(0)=P(t_j-), \hat{s}\in[0,1],\\
   \\
    \frac{d \xi_2^r}{d\hat{s}}=\gamma_r(\xi_2^r)l, \xi_2^r(0)=Q(t_j-), \hat{s}\in[0,1].
  \end{array} \right.
 \end{equation}

We construct the stochastic semi-implicit Euler scheme for $(\ref{3.1})$ which is shown as
 \begin{equation}\label{3.2}
\begin{split}
 P_{j+1}=P_j- \sigma_0(P_{j+1},Q_j) \Delta t_j-\sum_{r=1}^m \sigma_r(P_{j+1},Q_j)\diamond  \Delta L^r(t_j),\ P_0=P(t_0)=p,\\
 Q_{j+1}=Q_j+\gamma_0(P_{j+1},Q_j)\Delta t_j+\sum_{r=1}^m \gamma_r(P_{j+1},Q_j)\diamond  \Delta L^r(t_j),\ Q_0=Q(t_0)=q,
 \end{split}
\end{equation}
where $X_j:=(P_j,Q_j)$, $\Delta t_j=t_{j+1}-t_j$, $t_0<t_1<,...,<t_N$,   $\Delta L^r(t_j)=L^r(t_{j+1})-L^r(t_{j})$, and $j=0,1,...,N$. Here the Marcus integral for SDE(\ref{3.2}) is usually defined by
$$\sigma_r(P_{j+1},Q_j)\diamond  \Delta L^r(t_j)$$
$$=\sigma_r(P_{j+1},Q_j)\circ \Delta L_c^r(t_j)+\sigma_r(P_{j+1},Q_j)\Delta L_d^r(t_j)$$
$$+\sum_{t_j\leq s\leq t_{j+1}}\Big[\Phi_1^r(\Delta L^r(s),\sigma_r(P(s),Q(s)),P(s))-P(s)-\sigma_r(P(s),Q(s))\Delta L^r(s)\Big],$$
and
$$\gamma_r(P_{j+1},Q_j)\diamond  \Delta L^r(t_j)$$
$$=\gamma_r(P_{j+1},Q_j)\circ \Delta L_c^r(t_j)+\gamma_r(P_{j+1},Q_j)\Delta L_d^r(t_j)$$
$$+\sum_{t_j\leq s\leq t_{j+1}}\Big[\Phi_2^r(\Delta L^r(s),\gamma_r(P(s),Q(s)),Q(s))-Q(s)-\gamma_r(P(s),Q(s))\Delta L^r(s)\Big].$$
The notation $\circ$  denotes the Stratonovitch differential.
And the flow maps $\Phi_1^r(l,\sigma_r(P(s),Q(s)),P(s))$ and $\Phi_2^r(l,\gamma_r(P(s),Q(s)),Q(s))$ are the value at $s=1$ of the solutions defined through the ordinary differential equations, respectively,
\begin{equation}
 \left\{
  \begin{array}{lcr}
   \frac{d \xi_1^r}{ds}=\sigma_r(\xi_1^r)l, \xi_1^r(0)=P_{j+1}, s\in[0,1],\\
   \\
    \frac{d \xi_2^r}{ds}=\gamma_r(\xi_2^r)l, \xi_2^r(0)=Q_j, s\in[0,1].
  \end{array} \right.
 \end{equation}

We are in the position of the theorem which will show that the stochastic semi-implicit Euler scheme $(\ref{3.2})$ is symplectic, that is, the scheme $(\ref{3.2})$ preserves symplectic structure in the case of discrete time.

{\bf Theorem 3.2}\emph{
The scheme $(\ref{3.2})$ for the Hamiltonian SDE with multiplicative L\'{e}vy noise $(\ref{3.1})$ preserves symplectic structure.}

\begin{proof}

Due to the definition of symplectic structure, we only need to prove
$$dP_{j+1}\wedge dQ_{j+1}=dP_{j}\wedge dQ_{j},j=0,1,2,...,N,$$
where $dP_j$ and  $dQ_j$ are the differential of $P_j$ and $Q_j$, respectively.

 Take the differential with respect to $P$ of the first equation in SDE $(\ref{3.2})$, we obtain that
 $$ dP_{j+1}=dP_j- \frac{\partial \sigma_0}{\partial P}(P_{j+1},Q_j)\Delta t_j dP_{j+1}- \frac{\partial \sigma_0}{\partial Q}(P_{j+1},Q_j) \Delta t_jdQ_{j},$$
   $$ -\sum_{r=1}^m \frac{\partial \sigma_r}{\partial P}(P_{j+1},Q_j) dP_{j+1}\diamond  \Delta L^r(t_j) -\sum_{r=1}^m \frac{\partial \sigma_r}{\partial Q}(P_{j+1},Q_j) dQ_{j}\diamond  \Delta L^r(t_j),$$
 that is,
  $$\Big[\mathbb{I}+ \frac{\partial \sigma_0}{\partial P}(P_{j+1},Q_j)\Delta t_j+\sum_{r=1}^m \frac{\partial \sigma_r}{\partial P}(P_{j+1},Q_j)\diamond  \Delta L^r(t_j)\Big]dP_{j+1}$$
  $$=dP_j-\Big[\frac{\partial \sigma_0}{\partial Q}(P_{j+1},Q_j) \Delta t_j
   +\sum_{r=1}^m \frac{\partial \sigma_r}{\partial Q}(P_{j+1},Q_j) \diamond  \Delta L^r(t_j)\Big ] dQ_{j} $$
  where $\mathbb{I}$ is the $n \times n$ unit matrix.

  By the same way, we can obtain that
 $$dQ_{j+1}-\Big[\frac{\partial \gamma_0}{\partial P}(P_{j+1},Q_j) \Delta t_j+\sum_{r=1}^m \frac{\partial \gamma_r}{\partial P}(P_{j+1},Q_j)\diamond  \Delta L^r(t_j)\Big]dP_{j+1}$$
 $$=\Big [\mathbb{I}+ \frac{\partial \gamma_0}{\partial Q }(P_{j+1},Q_j)\Delta t_j + \sum_{r=1}^m \frac{\partial \gamma_r}{\partial Q}(P_{j+1},Q_j) \diamond  \Delta L^r(t_j)\Big ]dQ_{j}.$$
 Multiply the above two equations, we get
  $$\Big[\mathbb{I}+ \frac{\partial \sigma_0}{\partial P}(P_{j+1},Q_j)\Delta t_j+\sum_{r=1}^m \frac{\partial \sigma_r}{\partial P}(P_{j+1},Q_j)\diamond  \Delta L^r(t_j)\Big]dP_{j+1}\wedge dQ_{j+1}$$
  $$=\Big [\mathbb{I}+ \frac{\partial \gamma_0}{\partial Q }(P_{j+1},Q_j)\Delta t_j + \sum_{r=1}^m \frac{\partial \gamma_r}{\partial Q}(P_{j+1},Q_j) \diamond  \Delta L^r(t_j)\Big ]dP_{j}\wedge dQ_{j}.$$
  By the definition of Hamiltonian SDE $(\ref{3.1})$, we can obtain
  $$\frac{\partial \sigma_i}{\partial P}(P_{j+1},Q_j)= \frac{\partial \gamma_i}{\partial Q}(P_{j+1},Q_j)=\frac{\partial^2 H_i}{\partial Q \partial P}(P_{j+1},Q_j),i=0,1,...,m.$$

 As we know that the following inequality usually holds
$$\Big[\mathbb{I}+ \frac{\partial^2 H_0}{\partial Q \partial P}(P_{j+1},Q_j)\Delta t_j +\sum_{r=1}^m \frac{\partial^2 H_r}{\partial P \partial Q}(P_{j+1},Q_j)\diamond  \Delta L^r(t_j) \Big]\neq 0.$$

Therefore, we have
$$dP_{j+1}\wedge dQ_{j+1}=dP_{j}\wedge dQ_{j},j=0,1,2,...,N.$$
The proof of Theorem 3.2 is finished.
\end{proof}

\section{Convergence of symplectic Euler scheme}

{\bf Theorem 4.1}
\emph{If the inequality $0 < 1-18\sqrt{2}K^2(2\tau^2+12m^2\tau+m^2\eta)$ holds with Assumption 1, then the scheme $(\ref{3.2})$ for the Hamiltonian SDE with multiplicative L\'{e}vy noise $(\ref{3.1})$ based on one-step approximation is of the mean-square convergence order of accuracy 0.5, where $$\tau=\max_j\Delta t_j,  \eta=\max_j\eta_j,$$ and $\eta_j$ is the number of jumps in the time interval $[t_j,t_{j+1})$.}

\begin{proof}
We define
$$E_j=X_j-X_{t_j}.$$
According to the definition of the Marcus integral and $(\ref{3.3})$, we have
$$\mathbb{E}|E_{j+1}|^2=\mathbb{E}|X_{j+1}-{X}_{t_{j+1}}|^2
=\mathbb{E}\Bigg|
\begin{array}{c}
  P_{j+1}-P_{t_{j+1}}\\
  Q_{j+1}-Q_{t_{j+1}}
\end{array} \Bigg|^2\leq I_1+I_2+I_3,
$$
where
$$I_1=3\mathbb{E}\Bigg|
\begin{array}{c}
  P_j-P_{t_j}\\
  Q_j-Q_{t_j}
\end{array} \Bigg|^2,
I_2=3\mathbb{E}\Bigg|
\begin{array}{c}
\int_{t_j}^{t_{j+1}} [\sigma_0(P_{j+1},Q_j)-\sigma_0(P(s),Q(s))]ds \\
\\
\int_{t_j}^{t_{j+1}} [\gamma_0(P_{j+1},Q_j)-\gamma_0(P(s),Q(s))]ds
\end{array} \Bigg|^2
$$
and
$$I_3=3\mathbb{E}\Bigg|
\begin{array}{c}
\sum_{r=1}^m \int_{t_j}^{t_{j+1}}[\sigma_r(P_{j+1},Q_j)-\sigma_r(P(s),Q(s))] \diamond dL^r(s) \\
\\
\sum_{r=1}^m \int_{t_j}^{t_{j+1}}[\gamma_r(P_{j+1},Q_j)-\gamma_0(P(s),Q(s))] \diamond dL^r(s)
\end{array} \Bigg|^2.
$$
Now we need to estimate the value of $I_1$,$I_2$ and $I_3$,respectively.

Step 1.

To start, we have $I_1=\mathbb{E}|E_j|^2$.

Step 2.

Next, for $I_2$,we use the property of the norm and obtain
$$ \mathbb{E}\Big|\int_{t_j}^{t_{j+1}} [\sigma_0(P_{j+1},Q_j)-\sigma_0(P(s),Q(s))]ds \Big|^2\leq I_{21}+I_{22},$$
where
$$ I_{21}=2\mathbb{E}\Big|\int_{t_j}^{t_{j+1}} [\sigma_0(P_{j+1},Q_j)-\sigma_0(P(t_{j+1}),Q(t_{j+1}))]ds\Big |^2$$
and
$$I_{22}=2\mathbb{E}\Big|\int_{t_j}^{t_{j+1}} [\sigma_0(P(t_{j+1}),Q(t_{j+1}))-\sigma_0(P(s),Q(s))]ds\Big |^2.$$

Step 2(1).

On the one hand, we estimate $I_{21}$. It follows from the Cauchy-Scharz inequality that we have
$$ I_{21}\leq 2\Delta t_j\mathbb{E}\int_{t_j}^{t_{j+1}} \Big|\sigma_0(P_{j+1},Q_j)-\sigma_0(P(t_{j+1}),Q(t_{j+1}))\Big |^2ds$$
$$\leq 4\Delta t_j\mathbb{E}\int_{t_j}^{t_{j+1}} \Big[\Big|\sigma_0(P_{j+1},Q_j)-\sigma_0(P_{j+1},Q(t_{j+1}))\Big|^2$$
$$+\Big|\sigma_0(P_{j+1},Q(t_{j+1}))-\sigma_0(P(t_{j+1}),Q(t_{j+1}))\Big|^2 \Big]ds$$
$$\leq 4\Delta t_j K^2\mathbb{E}\int_{t_j}^{t_{j+1}}\Big[\Big|Q_j-Q(t_{j+1})\Big|^2+ \Big|P_{j+1}-P(t_{j+1})\Big|^2\Big]ds $$
$$\leq 4(\Delta t_j)^2 K^2\mathbb{E}\Big[2\Big|Q_j-Q_{j+1}\Big|^2+ 2\Big|Q_{j+1}-Q(t_{j+1})\Big|^2+\Big|P_{j+1}-P(t_{j+1})\Big|^2\Big]$$
\begin{equation}\label{4.0}
\leq 4(\Delta t_j)^2 K^2(C_1+3\mathbb{E}|E_{j+1}|^2),
\end{equation}
where $C_1$ depends on the assumption of the equations.

Step 2(2).

On the other hand, we estimate $I_{22}$. For $t\in [t_j,t_{j+1}]$ we obtain that
$$X_{t_{j+1}}=X_t- \int_{t}^{t_{j+1}} \sigma_0(X(s))ds-\sum_{r=1}^m \int_{t}^{t_{j+1}} \sigma_r(X(s)) \diamond  dL^r(s).$$
By the Chain rule of Marcus integral and the assumption that $\sigma_0$ is twice differentiable function, we have
$$\sigma_0(P_{t_{j+1}},Q_{t_{j+1}})-\sigma_0(P_t,Q_t)=- \int_{t}^{t_{j+1}} \sigma_0'(P(s),Q(s))\sigma_0(P(s),Q(s))ds$$
\begin{equation}\label{4.1}
-\sum_{r=1}^m \int_{t}^{t_{j+1}}\sigma_0'(P(s),Q(s))\sigma_r(P(s),Q(s)) \diamond  dL^r(s).
\end{equation}
According to Cauchy-Schwarz inequality, we get
$$I_{22}\leq 2 \Delta t_j \mathbb{E}\int_{t_j}^{t_{j+1}} \Big|\sigma_0(P(t_{j+1}),Q(t_{j+1}))-\sigma_0(P(s),Q(s))\Big |^2ds.$$
It is clear that from $(\ref{4.1})$ and Cauchy-Schwarz inequality we have
$$ \Big|\sigma_0(P(t_{j+1}),Q(t_{j+1}))-\sigma_0(P(s),Q(s))\Big |^2 \leq 2\Big[ \int_{t_j}^{t_{j+1}} \sigma_0'(P(s),Q(s))\sigma_0(P(s),Q(s))ds\Big]^2 $$
$$+2\Big [\sum_{r=1}^m \int_{t_j}^{t_{j+1}}\sigma_0'(P(s),Q(s))\sigma_r(P(s),Q(s)) \diamond  dL^r(s)\Big ]^2\leq C \Delta t_j,$$
where $C$ is a constant depends on $K$.
Then we can obtain that $$I_{22}\leq C\tau^2.$$
Therefore, we can have
\begin{equation}\label{4.2}
I_2\leq 3\sqrt{2}\tau^2\Big[4K^2(C_1+3\mathbb{E}|E_{j+1}|^2)+C\Big ].
\end{equation}

Step 3.

Last, for $I_3$, it follows from Cauchy-Scharz inequality that we have
$$\mathbb{E} \Big |\int_{t_j}^{t_{j+1}}[\sigma_r(P_{j+1},Q_j)-\sigma_r(P(s),Q(s))] \diamond dL^r(s)\Big |^2\leq I_{31}+I_{32},$$
where
$$ I_{31}=2\mathbb{E}\Big|\int_{t_j}^{t_{j+1}} [\sigma_r(P_{j+1},Q_j)-\sigma_r(P(t_{j+1}),Q(t_{j+1}))]\diamond dL^r(s)\Big |^2$$
and
$$I_{32}=2\mathbb{E}\Big|\int_{t_j}^{t_{j+1}} [\sigma_r(P(t_{j+1}),Q(t_{j+1}))-\sigma_r(P(s),Q(s))]\diamond dL^r(s)\Big |^2.$$

Step 3(1).

On the one hand, we estimate $I_{31}$. It follows from the definition of Marcus integral that we have
$$\mathbb{E}\Big |[\sigma_r(P_{j+1},Q_j)-\sigma_r(P(t_{j+1}),Q(t_{j+1}))]\diamond dL^r(s)\Big|^2\leq I_{311}+I_{312}+I_{313},$$
where
$$I_{311}=3\mathbb{E}\Big|\int_{t_j}^{t_{j+1}}\Big[\sigma_r(P_{j+1},Q_j)-\sigma_r(P(t_{j+1}),Q(t_{j+1}))\Big]\circ dL_c^r(s)\Big|^2,$$
$$I_{312}=3\mathbb{E}\Big|\int_{t_j}^{t_{j+1}}\Big[\sigma_r(P_{j+1},Q_j)-\sigma_r(P(t_{j+1}),Q(t_{j+1}))\Big]dL_d^r(s)\Big|^2$$
and
$$I_{313}=3\mathbb{E}\Big|\sum_{t_j\leq s\leq t_{j+1}}\Big[\Phi_1^r(\Delta L^r(s),\sigma_r(P_{j+1},Q_j),P_{j+1})$$
$$-\Phi_1^r(\Delta L^r(s),\sigma_r(P({t_{j+1}}),Q({t_{j+1}})),P({t_{j+1}}))\Big]\Big|^2.$$
It follows from Ito Isometry and $(\ref{4.0})$ that we have
$$I_{311}\leq 3\int_{t_j}^{t_{j+1}}\mathbb{E} \Big|\sigma_r(P_{j+1},Q_j)-\sigma_r(P(t_{j+1}),Q(t_{j+1}))\Big|^2 ds$$
 and
 $$I_{312}\leq 3\int_{t_j}^{t_{j+1}}\mathbb{E} \Big|\sigma_r(P_{j+1},Q_j)-\sigma_r(P(t_{j+1}),Q(t_{j+1}))\Big|^2 ds.$$
 Therefore, we have
 $$I_{311}+I_{312}\leq 6\int_{t_j}^{t_{j+1}}\mathbb{E} \Big|\sigma_r(P_{j+1},Q_j)-\sigma_r(P(t_{j+1}),Q(t_{j+1}))\Big|^2 ds $$
 $$\leq  12\tau K^2(C_1+3\mathbb{E}|E_{j+1}|^2).$$
By the Cauchy-Schwaz inequality and the Lipschitz condition, we have
$$I_{313}\leq 3 \eta_j\mathbb{E}\Big|\Phi_1^r(\Delta L^r(s),\sigma_r(P_{j+1},Q_j),P_{j+1})-\Phi_1^r(\Delta L^r(s),\sigma_r(P(t_{j+1}),Q(t_{j+1})),P(t_{j+1}))\Big|^2 $$
$$\leq 3\eta_jK^2\mathbb{E}|E_{j+1}|^2.$$
Therefore, for $I_{31}$, we have
$$I_{31}\leq 2(I_{311}+I_{312}+I_{313})$$
$$\leq 6K^2(\eta_j+12\tau)\mathbb{E}|E_{j+1}|^2+24\tau C_1K^2.$$

Step 3(2).

On the other hand, we estimate $I_{32}$.
$$I_{32}=2\mathbb{E}\Big|\int_{t_j}^{t_{j+1}} [\sigma_r(X(t_{j+1}))-\sigma_r(X(s))]\diamond dL^r(s)\Big |^2$$
$$\leq 2CK^2\mathbb{E} \Big | X(t_{j+1})-X(s)\Big |^2.$$
For $t\in [t_j,t_{j+1}]$ we obtain that
$$X_{t_{j+1}}=X_t- \int_{t}^{t_{j+1}} \sigma_0(X(s))ds-\sum_{r=1}^m \int_{t}^{t_{j+1}} \sigma_r(X(s)) \diamond  dL^r(s).$$
It is clear that we can obtain
$$\mathbb{E} \Big | X(t_{j+1})-X(s)\Big |^2\leq \Big[\int_{t_j}^{t_{j+1}} \sigma_0(X(s))ds-\sum_{r=1}^m \int_{t_j}^{t_{j+1}} \sigma_r(X(s)) \diamond  dL^r(s)\Big]^2\leq C\tau,$$
where $C$ is a constant depends on $K$.
Then we can obtain that $$ I_{32}\leq 2C^2K^2\tau.$$
Therefore, we can have
\begin{equation}\label{4.2}
I_3 \leq 3\sqrt{2}m^2\Big[6K^2(\eta_j+12\tau)\mathbb{E}|E_{j+1}|^2+24C_1K^2\tau+2C^2K^2\tau\Big].
\end{equation}

Step 4.

 All together, we have
 $$\mathbb{E}|E_{j+1}|^2\leq I_1+I_2+I_3.$$
That is,
$$\Big[1-18\sqrt{2}K^2(2\tau^2+12m^2\tau+m^2\eta)\Big]\mathbb{E}|E_{j+1}|^2\leq 3\mathbb{E}|E_{j}|^2+C\tau.$$

It follows from the assumption and the discrete version of Gronwall lemma that we have
$$\mathbb{E}|E_j|^2\leq C \tau.$$

Therefore, we have
$$\sup_{j\leq N}\|E_j\|=\sup_{j\leq N}\|X_j-X_{t_j}\|\leq C\tau^{\frac{1}{2}}.$$
The proof of Theorem 4.1 is finished.
\end{proof}

\section{Numerical implementation methods}
\subsection{Basic assumptions}
 
Let
\begin{equation}\label{4.00}
 L^r(t)=\sum_{k=1}^{\eta^r(t)}R_k^rH(t-\tau_k^r)+bW(t),r=1,2,...,m,
 \end{equation}
 where $\tau_k^r$ is the jump time with rate $\lambda$, $R_k^r\in \mathbb{R}$ is the jump size with distribution $\mu$, $\eta^r(t)$ is the number of jumps until time $t$, and
 $H(t)$ is the Heaviside function with unit jump at time zero.

 Due to the the realization of L\'{e}vy noises $(\ref{4.00})$, the Marcus integral for SDE(\ref{2.1}) through Marcus mapping is
 written as
$$X(t)=x+\int_0^tV_0(X(s))ds+b\sum_{r=1}^m\int_0^tV_r(X(s))\circ dW(s)$$
$$+\sum_{r=1}^{m}\sum_{k=1}^{\eta^r(t)}\Big[\Phi_g^r(X(\tau_k^r-),R_k^r)-X(\tau_k^r-)\Big],$$
where the flow map $\Phi_g^r$ at $t=\tau_k^r$ is defined through the ordinary differential equations
\begin{equation}\label{4.0}
 \left\{
  \begin{array}{lcr}
   \frac{d\xi^r}{ds}=V_r(\xi^r)R_k^r, s\in[0,1],\\
   \\
   \xi^r(0)=X(\tau_k^r-),\\
   \Phi_g^r(X(\tau_k^r-),R_k)=\xi^r(1).
  \end{array} \right.
 \end{equation}

Then we can simulate the orbits of Hamiltonian SDEs in a long time interval by symplectic Euler scheme, which will be shown in Section 5.2. Here we mainly consider the case $b=0$ in $(\ref{4.00})$, and we refer to the results which have been proposed in $\cite{T. Li, X. Wang}$ and $\cite{Q.Zhan6}$. And the realization of the case $b\neq 0$ will be considered in our future work.

\subsection{Pathwise symplectic Euler method}
We denote $\exp(\lambda)$ as the exponentially distributed random variable with mean $\frac{1}{\lambda}$. And we present this algorithm for the Hamiltonian SDE as follows,
\begin{equation}\label{4.1}
\begin{split}
 dP=-\sigma_0(P,Q)dt-\sum_{r=1}^m \sigma_r(P,Q)\diamond  dL^r(t),\ \ \  \ \ \ P(t_0)=p,\\
 dQ=\gamma_0(P,Q)dt+\sum_{r=1}^m \gamma_r(P,Q)\diamond  dL^r(t),\ \ \  \ \ \ Q(t_0)=q.
 \end{split}
\end{equation}

Step 1. Given $t=0$, initial value $(P_0,Q_0)$ and the end time $T$.

Step 2. Generate a waiting time $\tau\sim \exp(\lambda)$ and a jump size $R_r\sim \mu_r$, where $\mu_r(r=1,2,...,m)$ is the distribution of random jumps.

Step 3. Solve the following ODEs $(\ref{4.2})$ by symplectic Euler scheme with initial value $(P(t),Q(t))$ until time $s=\tau$ to get its solution $(P(u),Q(u)), u\in[t, t+\tau)$,
\begin{equation}\label{4.2}
\begin{split}
 dP=-\sigma_0(P,Q)dt-\sum_{r=1}^m \sigma_r(P,Q)\circ dW_t,\ \ \  \ \ \ P(t_0)=p,\\
 dQ=\gamma_0(P,Q)dt+\sum_{r=1}^m \gamma_r(P,Q)\circ dW_t,\ \ \  \ \ \ Q(t_0)=q.
 \end{split}
\end{equation}

Step 4. Solve the following ODEs $(\ref{4.3})$ with initial value $(P(t+\tau)-,Q(t+\tau)-)$ until time $s=1$ to get  $(P(t+\tau),Q(t+\tau))$,
\begin{equation}\label{4.3}
\begin{split}
&\frac{dx}{dt}=-\sum_{r=1}^m \sigma_r(P,Q)R_r,\ \ \  \ \ \ x(0)=P((t+\tau)-),\\
& \frac{dy}{dt}=\sum_{r=1}^m \gamma_r(P,Q)R_r,\ \ \  \ \ \ y(0)=Q((t+\tau)-).
 \end{split}
\end{equation}

Step 5. Set $t:=t+\tau$, and repeat Step 2 unless $t\geq T$.

\newpage

\section{ Numerical experiments}

We consider the following 2-dimensional SDE in the sense of Marcus $\cite{P. Wei}$,i.e., linear stochastic Kubo oscillator with multiplicative L\'{e}vy noise,
\begin{equation}\label{5.1}
\begin{split}
 & dP=-\alpha Qdt-\beta Q \diamond dL(t), \ P(t_0)=p ,\\
 & dQ=\alpha Pdt+\beta P \diamond dL(t), \ Q(t_0)=q,
 \end{split}
\end{equation}
where $\alpha$ and $\beta$ are constants, $L(t)$ is a one-dimensional L\'{e}vy noise, and $$H(P,Q)=p^2+q^2, H_1(P,Q)=p^2+q^2.$$
Obviously, it is a special nonlinear Hamiltonian SDEs with multiplicative L\'{e}vy noise which is the same as we discussed.
For any given initial values $(p,q)$, it follows from the results in $\cite{P. Wei}$ that the exact solution of SDE $(\ref{5.1})$ is
\begin{equation}\label{5.2}
\begin{split}
 & P(t)=p\cos(\alpha t+\beta L(t))-q\sin (\alpha t+\beta L(t)),\\
 & Q(t)=p\sin(\alpha t+\beta L(t))+q\cos(\alpha t+\beta L(t)).
 \end{split}
\end{equation}
The symplectic Euler scheme of SDE $(\ref{5.1})$ is written as
\begin{equation}\label{5.3}
\begin{split}
 & P_{j+1}=P_j-\alpha Q_j\Delta t_j-\beta Q_j \Delta L_j,\\
 & Q_{j+1}=Q_j+\alpha P_{j+1}\Delta t_j +\beta P_{j+1}\Delta L_j.
 \end{split}
\end{equation}

The explicit Euler scheme of SDE $(\ref{5.1})$ is written as
\begin{equation}\label{5.30}
\begin{split}
 & P_{j+1}=P_j-\alpha Q_j\Delta t_j-\beta Q_j \Delta L_j,\\
 & Q_{j+1}=Q_j+\alpha P_j\Delta t_j +\beta P_j\Delta L_j.
 \end{split}
\end{equation}

Section 6.1-6.2 are devoted to the preservation of symplectic structure and the convergence of the scheme $(\ref{5.3})$. And in the realization of the L\'{e}vy noise, we choose $L(t)$ to be a compound Poisson process with jump size which is simulated by the normal distribution $N(0,\sigma^2),\sigma=0.2$ and intensity $\lambda=5.0$.

\subsection{Preservation of symplectic structure of Hamiltonian SDE $(\ref{5.1})$ }

The results of our numerical experiments are shown as Fig.1-4, which includes three parts: the comparison of sample trajectories, the evolution of domains in the phase plane and the conservation of the Hamiltonian obtained by the scheme $(\ref{5.3})$, $(\ref{5.30})$ and the exact solution.

To start we apply the schemes $(\ref{5.3})$ and $(\ref{5.30})$  to Hamiltonian SDE $(\ref{5.1})$, and we can compare the oscillation of the numerical solutions obtained by the schemes $(\ref{5.3})$ and $(\ref{5.30})$ with the exact solutions from $(\ref{5.2})$. In order to improve the accuracy of the comparison, the initial conditions  are the same, that is, the step size is $dt=0.08$, $T=200.0$, $\alpha=0.1$, $\beta=0.1$, $N=2500.0$ and the initial values is $P(0)=0,Q(0)=1.0$. It is clear that the conditions of Theorem 4.1 are satisfied.
 \begin{figure}[h]
   \centering
   \begin{minipage}{6.5cm}
       \includegraphics[width=3.8in, height=2.80in]{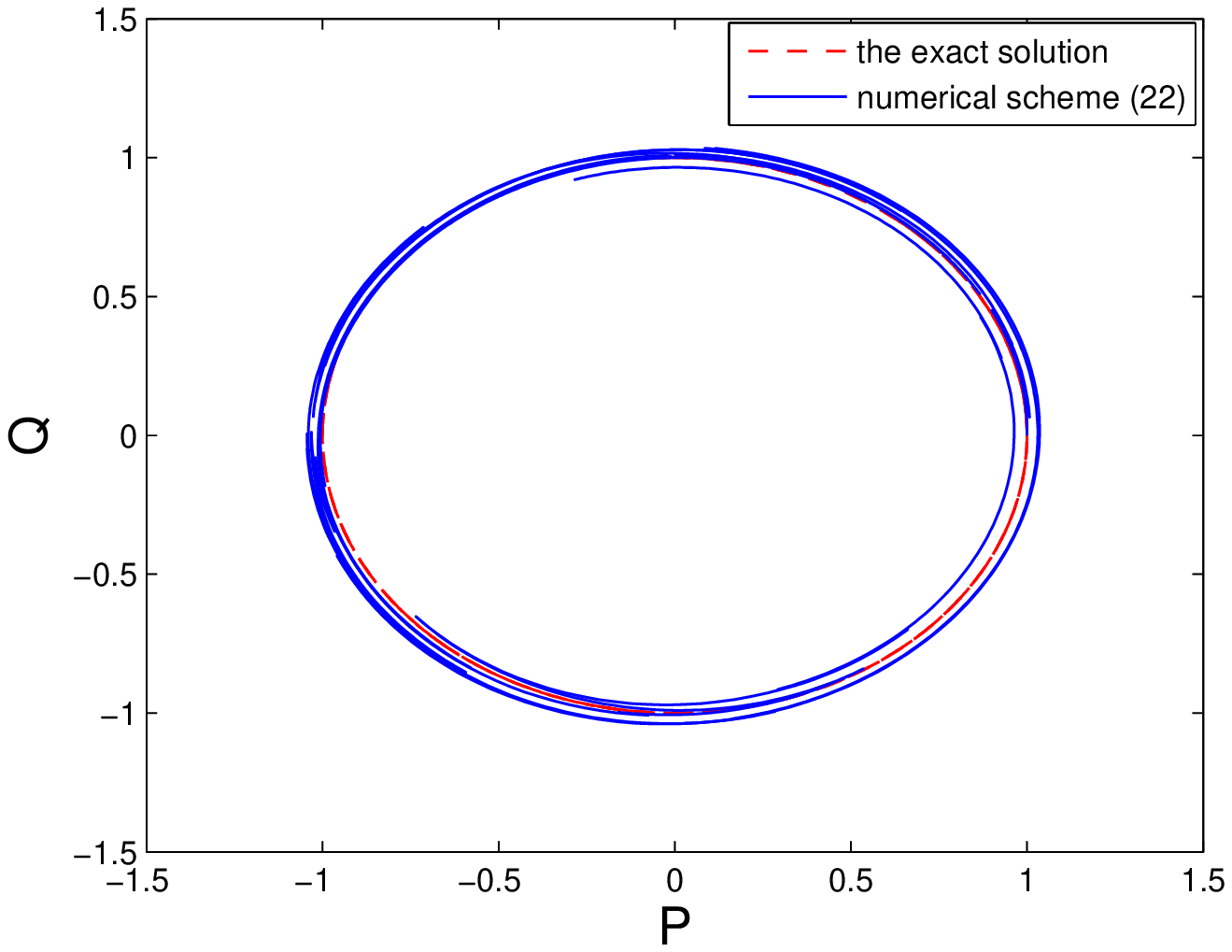}
    \end{minipage}
    \begin{minipage}{6.5cm}
       \includegraphics[width=3.8in, height=2.80in]{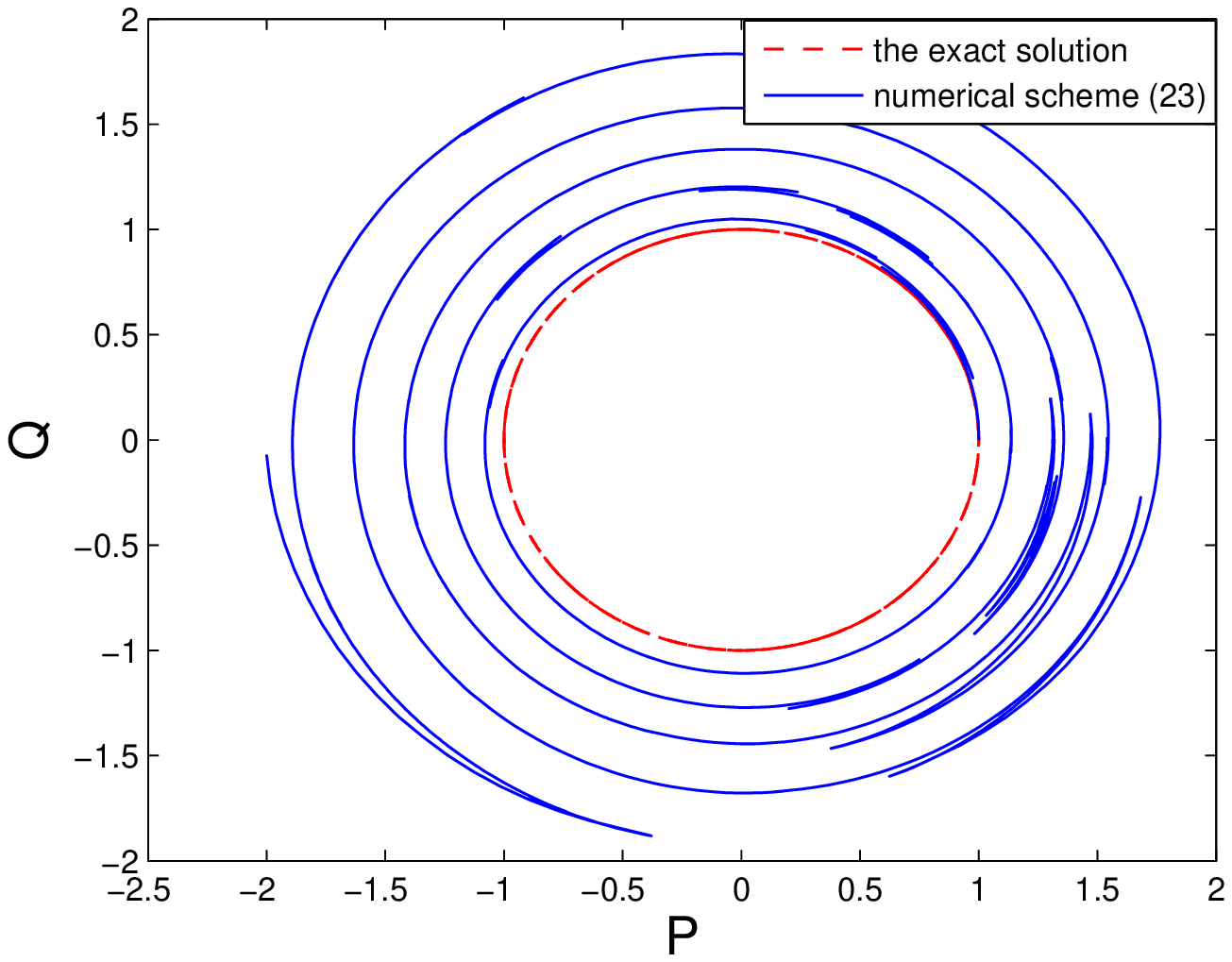}
    \end{minipage}
    \par{\scriptsize  Fig.1. Comparison of zoom in parts of the exact orbit of the solutions to SDE $(\ref{5.1})$ obtained by $(\ref{5.2})$  and a sample orbit obtained by the scheme $(\ref{5.3})$(Left) and $(\ref{5.30})$ (Right), respectively.}
\end{figure}

As we can see from Fig.1 that the approximations of a sample orbit of Hamiltonian SDE $(\ref{5.1})$ are simulated by the symplectic method, the scheme $(\ref{5.3})$, as well as the non-symplectic method, the scheme $(\ref{5.30})$, respectively. The exact phase trajectory $(\ref{5.2})$ is obtained, too.

  It shows the fact that in the time interval $[0,200.0]$, the orbit of the exact solution coincides almost well with that obtained by the scheme $(\ref{5.3})$, which is demonstrated in the left panel of Fig.1, while the orbit obtained by the scheme $(\ref{5.30})$ does not circle that of the exact solution, it disperses spirally  and quickly from the latter, which is shown as the right panel of Fig.1. It is obvious that the scheme $(\ref{5.3})$ has higher performance to preserve the circular orbits than the scheme $(\ref{5.30})$. That is,  the structure of the orbit of the solution to SDE $(\ref{5.1})$ obtained by the scheme $(\ref{5.30})$  obviously does not conserve the circular structure of that of the exact solution.
 \begin{figure}[htbp]
\centering
\includegraphics[width=3.8in, height=2.60in]{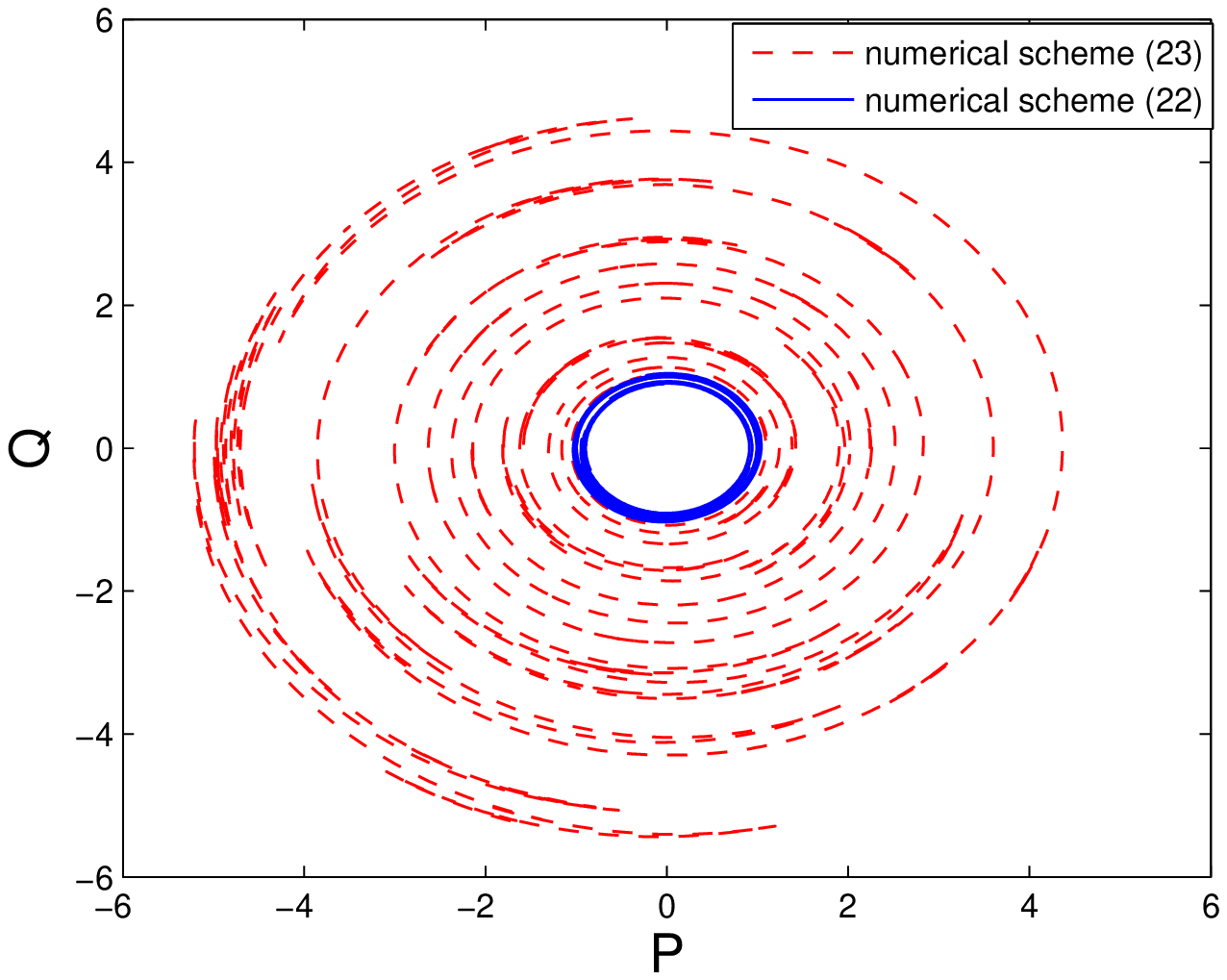}
\par{\scriptsize  Fig.2. Comparison of a sample trajectory of SDE $(\ref{5.1})$ obtained by symplectic method, the scheme $(\ref{5.3})$(blue line) and by the non-symplectic method,the scheme $(\ref{5.30})$(red line), respectively..}
\end{figure}

 \begin{figure}[h]
   \centering
   \begin{minipage}{6.5cm}
       \includegraphics[width=3.8in, height=2.60in]{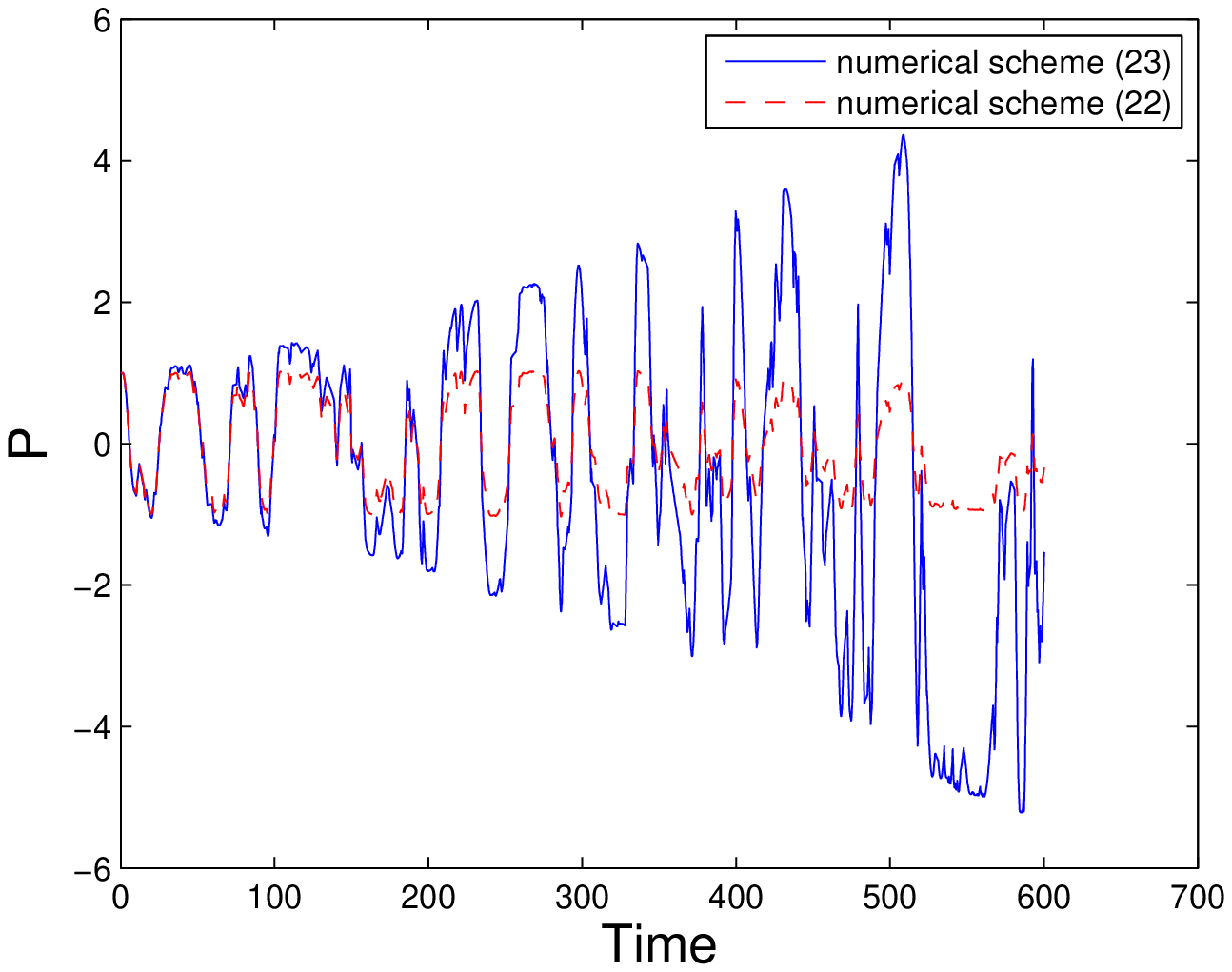}
    \end{minipage}
    \begin{minipage}{6.5cm}
       \includegraphics[width=3.8in, height=2.60in]{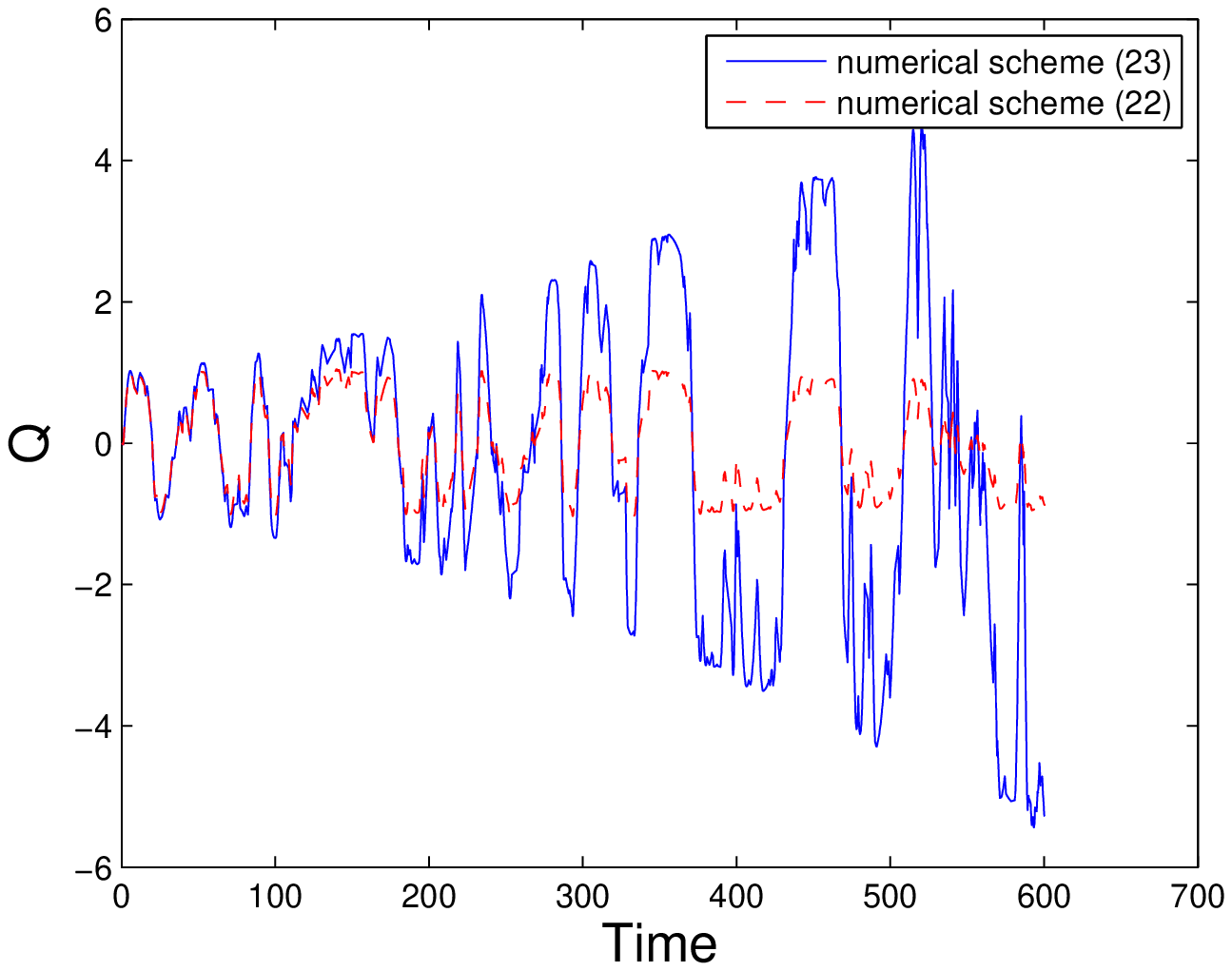}
    \end{minipage}
    \par{\scriptsize   Fig.3. Comparison of the oscillation in the directions $P$ and $Q$ of a sample trajectory of SDE $(\ref{5.1})$ obtained by symplectic method, the scheme $(\ref{5.3})$(red line) and by the non-symplectic method,the scheme $(\ref{5.30})$(blue line), respectivelyy.}
\end{figure}

 As we can see from Fig. 2 and Fig. 3, the oscillation in the directions $P$ and $Q$ of the trajectory of numerical solution obtained by the scheme $(\ref{5.3})$ is much smaller than that by the scheme $(\ref{5.30})$. And it is obvious that it disperses spirally  and quickly from the latter in the time interval $[0,600.0]$.

These results indicate that the scheme $(\ref{5.30})$ is unsuitable to simulate Hamiltonian SDE $(\ref{5.1})$ in a long time interval. In contrast to the scheme $(\ref{5.30})$, the scheme $(\ref{5.3})$ reproduces the trajectory of SDE $(\ref{5.1})$ more accurately.

Next we check the Hamiltonians of SDE $(\ref{5.1})$.
 \begin{figure}[htbp]
   \centering
   \begin{minipage}{6.5cm}
       \includegraphics[width=3.8in, height=2.60in]{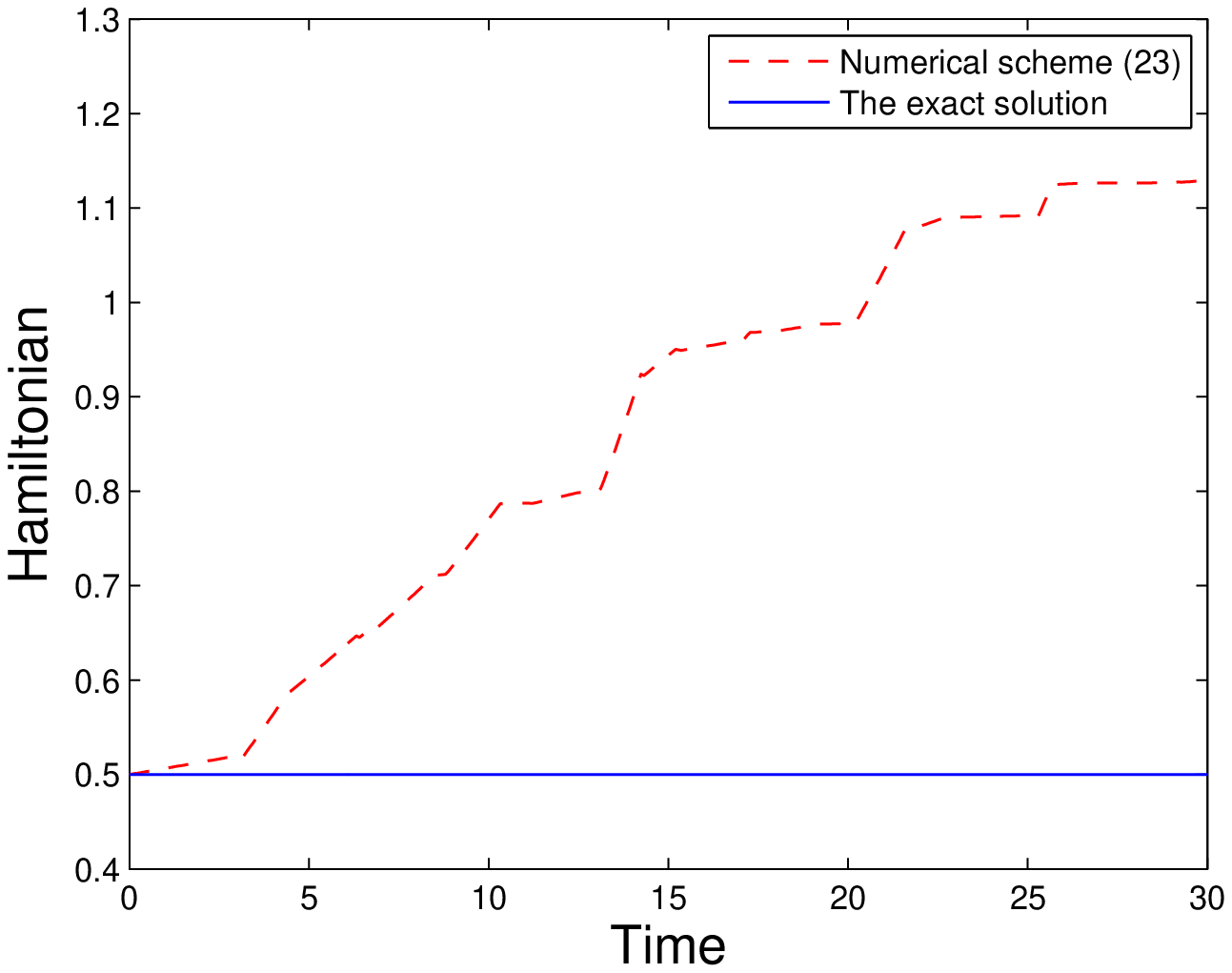}
    \end{minipage}
    \begin{minipage}{6.5cm}
       \includegraphics[width=3.8in, height=2.60in]{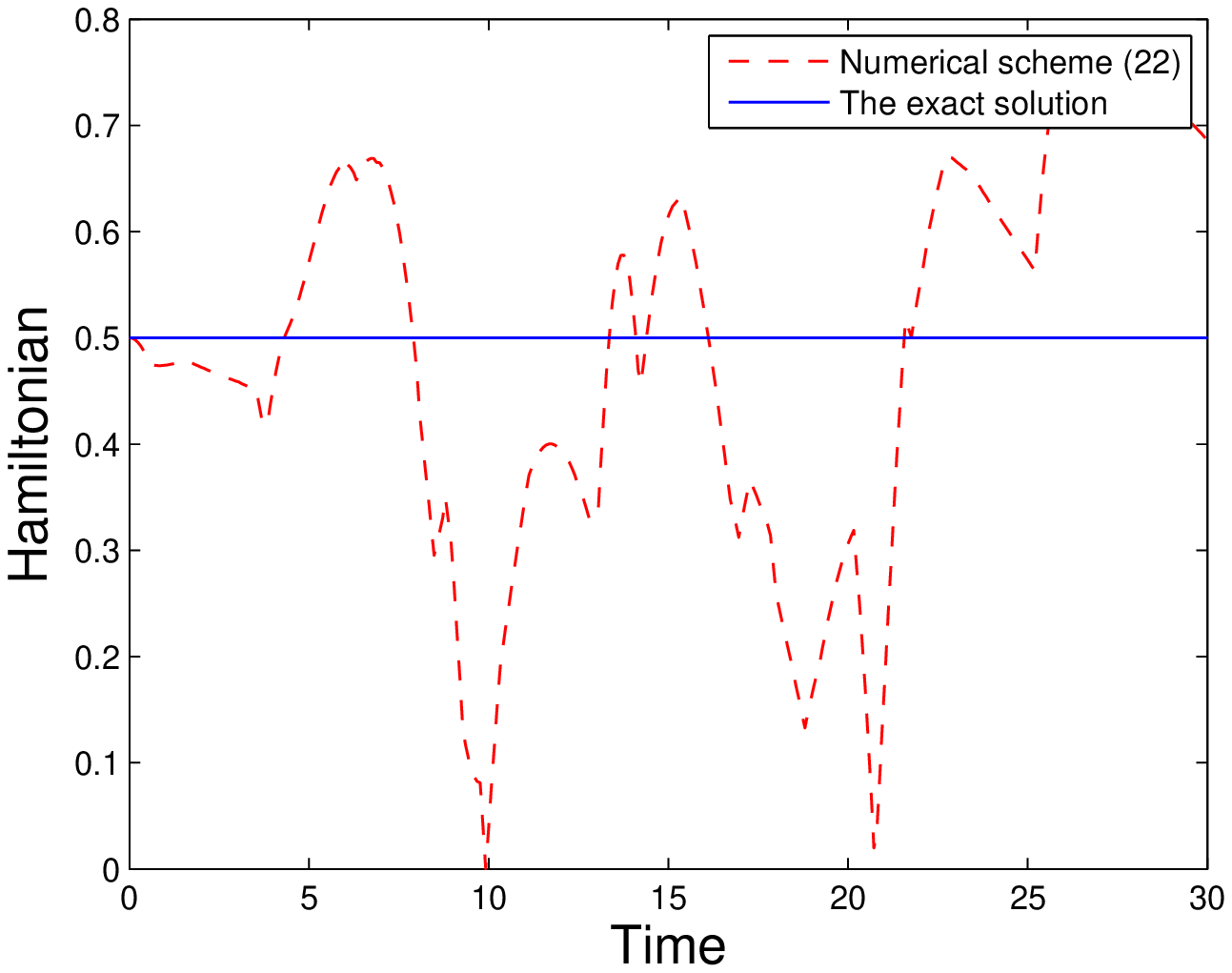}
    \end{minipage}
    \par{\scriptsize   Fig.4.  Conservation of the Hamiltonian of numerical solution by the scheme $(\ref{5.3})$ and $(\ref{5.30})$, respectively.}
\end{figure}

It can be seen from Fig.4 that the Hamiltonian, $H(P,Q)$, is an invariant of the exact solution of SDE $(\ref{5.1})$. As it shows that the curve of Hamiltonian jumps around the line Hamiltonian$=0.5$, which demonstrates it can be approximately preserve by the scheme $(\ref{5.3})$ because of the L\'{e}vy noise. Here approximate preservation of the Hamiltonian means that there is a bounded oscillation around the Hamiltonian of the exact solution in the discrete time case.   However, non-symplectic numerical scheme, the scheme $(\ref{5.30})$ dose not has this property such that the Hamiltonian increases indefinitely, which is illustrated as Fig.4.

\subsection{Convergence of the scheme $(\ref{5.3})$}

This numerical experiment examines the convergence of the scheme $(\ref{5.3})$. It is not difficult to see from Fig.5 that the convergence rate satisfies the inequality  log($\|$error $\|$)$\leq 0.5$ for the end time $T=100.0$. Because of the restrictions on computation capability, we only consider the case in the interval $[0,100]$. Due to discontinuous inputting of the L\'{e}vy noise, the curve has some jumps in some uncertain time moments, but it almost lays down the straight line log($\|$error$\|$)$=0.5$. And  These phenomena verify the results of Theorem 4.1 that the mean-square order of the proposed method is 0.5. In this test we choose the same parameters as Section 5.1, the mean-square norm is taken as $(\ref{2.3})$.
 \begin{figure}[htbp]
\centering
\includegraphics[width=3.8in, height=2.60in]{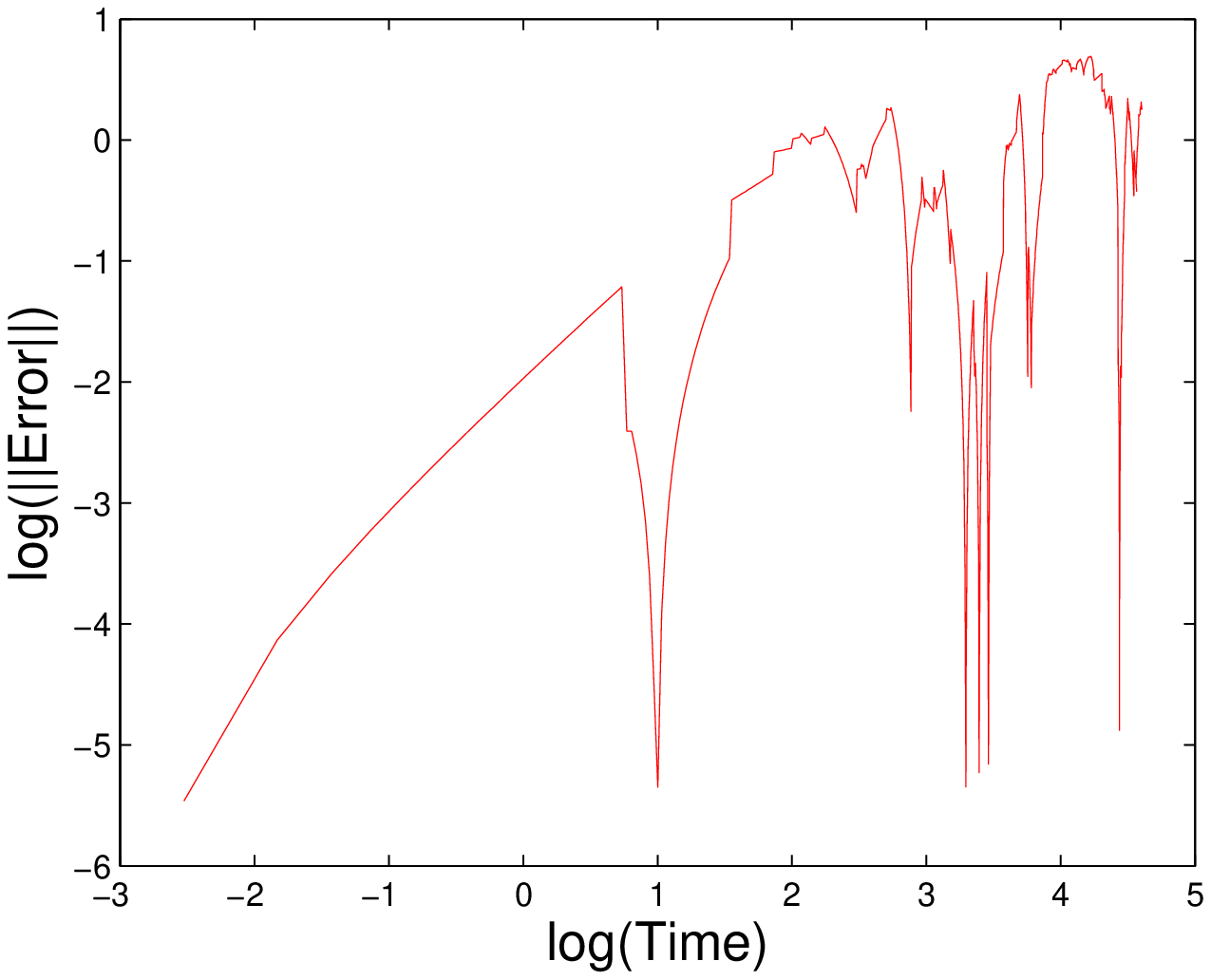}
\par{\scriptsize  Fig.5. The mean-square convergence rate of the scheme $(\ref{5.3})$.}
\end{figure}

\section{Conclusion}
 This paper mainly focuses on the construction, the convergence analysis and the numerical implementation of symplectic Euler scheme for Hamiltonian SDEs with multiplicative L\'{e}vy noise in the Marcus form. Much attention are paid to the numerical realization of this symplectic Euler scheme. We show that the main results and numerical implementation methods can be applied to the numerical simulations of the dynamical behaviour of Hamiltonian SDEs with multiplicative L\'{e}vy noise in the sense of Marcus. The results show that the method is effective and the numerical experiments are performed and match the results of theoretical analysis almost perfectly by comparing with non-symplectic method.

\section*{Statements}
All data in this manuscript is available. And all programs will be available on the WEB GitHub\cite{Q.Zhan7}.

\section*{Acknowledgments}

This work is supported by NSFC(No. 61841302 and 11771449). This work is also supported by the Science Research Projection in the Education Department of Fujian Province, No. JT180122, Education Reform Fund of Fujian Agriculture and Forestry University, No. 111418136. Qingyi Zhan would like to thank the Department of Applied Mathematics, Illinois Institute of Technology for the hospitality during his visit(2019-2020). And he wants to
thank Prof. Jinqiao Duan, Prof. Xiaofan Li and all members in Lab. for stochastic dynamics and computation in IIT for many fruitful discussions
during that period. Qingyi Zhan would also like to acknowledge
the sponsorship of the China Scholarship Council, CSC No. 201907870004.

\begin{backmatter}

\section*{Competing interests}
The authors declare that they have no competing interests.

\section*{Author details}
1.College  of Computer and Information Science, Fujian Agriculture and Forestry  University, Fuzhou, Fujian, 350002,China;

2.Department of Applied Mathematics, Illinois Institute of Technology, Chicago, IL, 60616, USA ; \\

3.College of Hydropower and Information Engineering, Huazhong University of Science and Technology, Wuhan,430074,China;\\
*:Corresponding Author:Q.Zhan: zhan2017@fafu.edu.cn. Co-author: Jinqiao Duan:duan@iit.edu; Xiaofan Li: lix@iit.edu;Yuhong Li: liyuhong@hust.edu.cn.


\bibliographystyle{bmc-mathphys} 
\bibliography{bmc_article}      


\end{backmatter}
\end{document}